\documentclass[a4paper, 11pt]{article}
\usepackage{amsmath,amsthm,amsfonts,amssymb}
\theoremstyle{plain}

\numberwithin{equation}{section}

\newtheorem{theorem}[equation]{Theorem}

\newtheorem*{teorema}{Theorem A}
\newtheorem*{teoremab}{Theorem B}
\newtheorem{proposition}[equation]{Proposition}

\newtheorem{lemma}[equation]{Lemma}

\newtheorem{example}[equation]{Example}

\newcommand{\Irr}{\operatorname{Irr}}

\newcommand{\Ker}{\operatorname{Ker}}

\newcommand{\dl}{\operatorname{dl}}

\theoremstyle{definition}

\newtheorem{definition}[equation]{Definition}

	\title{Products of characters and derived length II}

\author{}

\date{}

\begin{document}
\maketitle
\centerline{by}

\centerline{Edith Adan-Bante}
\centerline{Department of Mathematics}
\centerline{University of Illinois at Urbana-Champaign}
\centerline{Urbana, IL 61801 USA}

Currently at the University of Southern Mississippi Gulf Coast,
730 East Beach Blvd, Long Beach, MS 39560

\noindent

The author was partially supported by the 
National Science Foundation by
grant DMS 99-70030.

\newpage

\begin{section}{Introduction}

	Let $G$ be a finite group. Let $\chi$ and $\psi$
 be irreducible
complex characters, i.e. 
 irreducible characters over
the complex numbers. Since a product of characters is a character,
$\chi \psi$ is a character.
 Then the decomposition of
the character $\chi \psi$ into its distinct 
irreducible 
constituents $ \alpha_1, \ \alpha_2, \ldots, \alpha_n$
has the form
\begin{equation*}	
 \chi\psi=  \sum_{i=1}^n a_i \alpha_i
\end{equation*}
\noindent where $n>0$ and 
$a_i>0$ is the multiplicity of 
$\alpha_i$ for each $i=1, \ldots, n$. 
Let  $\eta(\chi\psi)= n $ be the
number of distinct 
 irreducible constituents 
of  the character $\chi \psi$. Denote by $\Ker(\chi\psi)$ the kernel of the 
character $\chi\psi$. Also denote by $\dl(\Ker(\alpha)/\Ker(\chi \psi))$ 
the derived
length of the group $\Ker(\alpha)/\Ker(\chi \psi)$.

 Given $\chi\in \Irr(G)$, let $\overline{\chi}$ 
 be the complex conjugate of $\chi$,
  i.e. $\overline{\chi}(g)=\overline{\chi(g)}$ for all $g \in G$. 
 In Theorem A of \cite{edith}, it is proved that there exist universal
constants
 $C_0$ and $D_0$ such that for any  finite solvable group $G$
 and any $\chi\in \Irr(G)$ we have that 
 $\dl(G/\Ker(\chi))\leq C_0\eta(\chi\overline{\chi})+D_0$.
The following generalizes that result and it is the main result of 
this paper.

\begin{teorema}
 There exist universal constants $C$ and $D$ such that for
any finite solvable group $G$, any irreducible 
characters $\chi$ and $\psi$ 
of $G$, and  any   irreducible
constituent   $\alpha$         of $\chi \psi$, we have 
\begin{equation*}
\dl(\Ker(\alpha)/\Ker(\chi \psi)) \leq  C\eta(\chi\psi ) + D.
\end{equation*}
\end{teorema}

In Theorem \ref{supersolvable} we prove that
we may take 
$C=2$ and $D=-1$ 
if the group $G$ in Theorem A is,
 in addition, supersolvable.

\begin{teoremab}
Let $G$ be a finite group 
and $\chi,\psi\in \Irr(G)$ be characters.
If $(\chi(1),\psi(1))=1$, then $\dl(\Ker(\alpha)/\Ker(\chi\psi))\leq 1$
for any irreducible constituent $\alpha$ of the product $\chi\psi$.
\end{teoremab}
 
Thus if $(\chi(1), \psi(1))=1$, then the irreducible
constituents of the product $\chi\psi$ are ``almost''
faithful characters of the group $G/\Ker(\chi\psi)$. 

{\bf Acknowledgment.} I  thank Professor Everett C. Dade 
for helpful discussions. I also thank the referee for
 her/his suggestions for improving the 
presentation of this work, and especially for shortening 
the proof of Theorem A.
 
\end{section}
\begin{section}{Proof of Theorem A}

We will be using the notation of \cite{isaacs}. 
In addition, we write $M\subseteq N$ if $M$ is a subgroup of $N$. Also
$M\subset N$ denotes that $M$ is a proper subgroup of $N$. 

\begin{definition}
	Assume that $G$ acts on a finite vector space $V$. 
We define  $m(G,V)$ as the number of  orbits of nonzero vectors
under the action of $G$.
\end{definition}

	The key tool for the proof of Theorem A is the following.
\begin{lemma} [Keller] \label{keller}
There exist universal constants $A$ and $B$ such that for 
any finite solvable group $G$ acting faithfully  on
a finite vector space $V$ we have that
$$\dl(G)\leq A m(G,V) +B.$$
\end{lemma}
\begin{proof} 
See \cite{keller}.
\end{proof}

We  emphasize  that the result of Keller is much stronger than
Lemma \ref{keller}. He proved that the derived length of $G$
is bounded by a logarithmic function of the number of different
sizes of $G$-orbits on $V$.

\begin{definition}\label{orbits} We define the function $f$ by
\begin{equation*}
 f(n)= A n + B,
\end{equation*}
\noindent for any integer $n\geq 1$, 
where $A$ and $B$ are as in Lemma \ref{keller}.
\end{definition}

\begin{definition}\label{extreme}
Let $G$ be a finite group and $M$ and $N$ be normal 
subgroups of $G$ such that $M\subset N$. Let $\theta\in \Irr(N)$. We say that 
$(N,M, \theta)$ is an {\bf extreme triple} of $G$ if 
for any normal subgroup $K$ of $G$ such that
$M\subset K\subseteq N$, we have that $\theta_K$ is irreducible but 
$\theta_M$ is reducible.
\end{definition}
Observe that given $N\trianglelefteq G$ and 
 a nonlinear character $\theta\in \Irr(N)$ of $N$, we can always find
$M\triangleleft
G$ such that $(N,M,\theta)$ is an extreme triple.

\begin{lemma}\label{referea} 
Let $(N,M, \theta) $ be an extreme triple of 
$G$ and suppose $M\subset L\subseteq N$,
where $L\triangleleft G$ and $L/M$ is abelian. Then
${\bf C}_{N/M}(L/M)$ is abelian.
\end{lemma}
\begin{proof}
By definition of extreme triple, we have that $\theta_L\in \Irr(L)$
and $\theta_M$ is reducible. Since $L/M$ is abelian, by Theorem
6.22 of \cite{isaacs} we have that there exists a subgroup $U$ 
containing $M$ such that $U$ has prime index in $L$ and $\theta_U$
is reducible. Since $\theta_U$ is reducible and $|L:U|=p$, we have that
$\theta_U$ is a sum of $p$ distinct irreducible characters of $U$. Let
$\varphi\in \Irr(U)$ be an irreducible constituent of $\theta_U$. 

Let $C={\bf C}_N(L/M)$. Note that $L\subseteq C$ and thus $U\triangleleft C$
and
 $\theta_C\in \Irr(C)$. It follows that $|C: C_{\varphi}|=p$ and
 $C=LC_{\varphi}$. Since $L$ centralizes $C/M$, we see that 
$C_{\varphi}\trianglelefteq C$ and  
$C/C_{\varphi}$ is abelian. Then $MC'\subseteq C_{\varphi}$ and
thus $\theta_{MC'}$ is reducible. Since 
$MC'\trianglelefteq G$, it follows that $MC'=M$. Thus
 $C'\subseteq M$ as wanted.
\end{proof} 

We introduce more notation. If $N/M$ is a normal section of $G$
and  $\Delta$ is a character of $G$, we write  $S_{\Delta} (N/M)$
to denote the (possibly empty) set of those irreducible 
constituents $\alpha$ of $\Delta$ such that $M\subseteq \Ker(\alpha)$ and
$N\not\subseteq \Ker(\alpha)$.

\begin{lemma}\label{refereb}
Let $(N,M,\theta)$ be an extreme triple of $G$, where $N/M$ is 
solvable. Assume that $\theta\overline{\theta}$ is 
(not necessarily irreducible) constituent of $\Delta_N$, for
some character
$\Delta$ of $G$. Then $S_{\Delta} (N/M)$ is nonempty, and
the derived length of $N/M$ is at most $1+f(|S_{\Delta}(N/M)|)$, 
where $f$ is the function as in Definition  \ref{orbits}.
\end{lemma}
\begin{proof}
Since $N/M$ is solvable and non-trivial, we can choose $L\trianglelefteq G$
with $M\subset L\subseteq N$ such that $L/M$ is a chief factor of $G$.
 Also, as  in the 
previous lemma, we know that there exists a subgroup $U$ containing 
$M$ and of prime index $p$ in $L$
such that $\theta_U$ is a sum of $p$ distinct  
irreducible constituents.
Thus $\theta$ vanishes on $L\setminus U$. Since $\theta_L$ is $N$-invariant,
$\theta $ vanishes on $L\setminus D$, where $D$ is the intersection of 
the $N$-conjugates of $U$. Since $M\subseteq D\subset L $ and $D\triangleleft
N$,
we can choose a chief factor $L/K$ of $N$ such that $D\subseteq K$,
and thus $\theta$ vanishes on $L\setminus K$.

Write $\varphi=\theta_L\in \Irr(L)$ and 
let $\lambda$ be a nonprincipal linear character of $L/K$. Then 
$\varphi\lambda=\varphi$, and this $\lambda$ is a constituent of
$\varphi\overline{\varphi}$, which in turn is a constituent 
of $\Delta_L$. It follows that $\lambda$ is a constituent of
$\alpha_L$, where $\alpha$ is an irreducible constituent of
$\Delta$. In particular, we see that $\alpha\in  S_{\Delta} (N/M)$,
which is therefore nonempty.

Set $H={\bf N}_G(K)$. Observe that  $H\supseteq N$. 
All members of the $H$-orbit of $\lambda$ are linear constituents
of $\alpha_L$ having kernels containing $K$. Conversely, 
we show that if 
$\nu$ is any linear constituent of $\alpha_L$ such that $K\subseteq \Ker(\nu)$,
then $\nu$ lies in the $H$-orbit of $\lambda$. 
Certainly $\nu= \lambda^g$ for some element $g\in G$, and 
thus $\nu$ is nonprincipal. Also, since $K\subseteq \Ker(\lambda)$,
we see that $K^g\subseteq \Ker(\nu)$. We assumed that
$K \subseteq \Ker(\nu)$, and thus $K K^g\subseteq \Ker(\nu)\subset L$.
Observe that $K K^g\triangleleft N$ since $K^g\subseteq N^g=N$ and $K$ is normal
in $N$.
Since $KK^g\triangleleft N$ and $L/K$ is a chief factor of $N$,
it follows that $K=K^g$ and thus $g\in H$.
We see now that for each $H$-orbit of nonprincipal linear characters of 
$L/K$, there is a member of $S_{\Delta} (N/M)$
that lies over all members of this orbit and over no other linear 
character of $L/K$. It follows that the number of $H$-orbits on $L/K$,
which is the same as the number of $H$-orbits of linear 
characters of $L/K$, is at most
$|S_{\Delta} (N/M)|$. By Lemma \ref{keller}, therefore, the derived
length of
$H/{\bf C}_H(L/K)$ is at most $f(|S_{\Delta} (N/M)|)$. 
Set $m=f(|S_{\Delta} (N/M)|)$. Then the derived length of $N/{\bf C}_N(L/K)$
is at most $m$ and thus $N^{(m)}\subseteq {\bf C}_N(L/K) $. 
 Since $[L,N^{(m)}]\subseteq K$, 
 $K\subset L$, $L/M$ is a chief factor of $G$ and 
 $[L,N^{(m)}]\triangleleft G $, we see that  $[L,N^{(m)}]\subseteq M$. 
Thus  $N^{(m)}$ centralizes $L/M$.
Since ${\bf C}_{N/M}(L/M)$ is abelian by the previous lemma, 
it follows that $N^{(m+1)}\subseteq M$, as desired.
\end{proof}
\begin{lemma}\label{partb}
Let $G$ be a finite group and $\chi, \psi \in \Irr(G)$. 
Let $\alpha \in \Irr(G)$ be an irreducible constituent 
of the product $\chi \psi$ and $K= \Ker(\alpha)$. 
Then there exists some 
 $\theta \in \Irr(K)$ such that 
$[\chi_K, \theta]\neq 0$ and $[\overline{\psi}_K, {\theta}]\neq 0$.
\end{lemma}
\begin{proof}
Since $K= \Ker(\alpha)$, we have that $\alpha_K= \alpha(1) 1_K$. Thus
$$\alpha(1)[(\chi\psi)_K, 1_K]= [(\chi\psi)_K, (\alpha)_K] >0.$$
Therefore there exist $\,\theta, \sigma\in \Irr(K)$ 
 such that 
$[\chi_K, \theta]\neq 0$, $[\psi_K, \sigma]\neq 0$ and 
$[ \theta\sigma, 1_K] \neq 0$. Since
$\theta$ and $\sigma$ are irreducible characters and 
 $[ \theta\sigma , 1_K] = [\theta, \overline{\sigma}]>0$, we conclude that
$\theta = \overline{\sigma}$. Thus  $[\overline{\psi}_K, \theta]\neq 0$.    
\end{proof}
\begin{proof}[Proof of Theorem A]
Write $\Delta=\chi\psi$ and $n=\eta(\chi\psi)$.
Note that we can assume that $\Delta$ is
faithful. Let $K=\Ker(\alpha)$ so that our task is to show that the
 derived length $\dl(K)$ is bounded by a linear function on $n$.
By Lemma \ref{partb}, we can choose a character $\nu \in \Irr(K)$ that lies 
under both 
$\chi$ and $\overline{\psi}$.  Also we see that if $N\triangleleft G$ 
and $N\subseteq \Ker(\nu)$, then $N\subseteq \Ker(\Delta)$ and thus 
$N=1$. 

Write $K_0=K$. If $\nu$ is nonlinear, we can choose a subgroup $L$ 
such that $(K,L,\nu)$ is an extreme triple of $G$. Write $K_1=L$ and 
let $\nu_1$ be an irreducible constituent of $\nu_L$. If
$\nu_1$ is nonlinear, we can repeat the process and choose an extreme
triple of 
$G$ of the form $(K_1, K_2, \nu_1)$. Continuing like this, we 
obtain a series $K=K_0\supset K_1\supset \cdots\supset K_r$ of 
normal subgroups of $G$ and characters $\nu_i\in \Irr(K_i)$ such 
that $(K_i, K_{i+1}, \nu_i)$ is an extreme triple  for $0\leq i<r$ 
and where $\nu_r$ is linear. Also, we know that for each subscript 
$i$, the 
character $\nu_i$ is a constituent of $\nu_{K_i}$, and thus $\nu_i$ 
lies under both $\chi$ and $\overline{\psi}$. It follows that 
$\nu_i \overline{\nu_i}$ is a constituent of $\Delta_{K_i}$.  
Therefore by Lemma \ref{refereb}, the sets $S_{\Delta}(K_i/K_{i+1})$ 
are nonempty for $0\leq i<r$. Since they are certainly 
disjoint and $\alpha$ lies in none of the sets
$S_{\Delta}(K_i/K_{i+1})$,
 we have that $r<n$. Also, writing 
$s_i=|S_{\Delta}(K_i/K_{i+1})|$, we see that $\sum s_i< n$.

By Lemma \ref{refereb}, the derived length of $K_i/K_{i+1}$ is at most 
$1+f(s_i)$. Also, because $\nu_r$ is linear, we see that 
$(K_r)'\subseteq \Ker(\nu_r)$. Since $(K_r)'\triangleleft G$,
it follows that $(K_r)'\subseteq \Ker(\Delta)$, and hence
$(K_r)'=1$ and $K_r$ is abelian. We see now that the derived
length $\dl(K)$ is at most $1+\sum_{i=0}^{r-1}(1+f(s_i))$.
Since $f(x)=Ax+B$, we have that $\dl(K)\leq A\sum s_i +rB +r +1$.
But $\sum s_i< n$ and $r<n$, and so $\dl(K)< (A+B+1)n+1$.
The proof is now complete.
\end{proof}

\begin{theorem}\label{supersolvable}
Let $G$ be a finite supersolvable group and $\chi\psi\in \Irr(G)$
be characters of $G$. Let $\alpha \in \Irr(G)$ be any  
constituent  of the product $\chi \psi$. Then  
\begin{equation*}
\dl(\Ker(\alpha)/\Ker(\chi \psi)) \leq  2\eta(\chi \psi ) - 1.
\end{equation*}
\end{theorem}
\begin{proof} Set $K=\Ker(\alpha)$ and $n=\eta(\chi \psi)$.
 We choose a character $\nu\in \Irr(K)$ that lies
under both $\chi$ and $\overline{\psi}$. 
As in the proof of Theorem A, 
 we 
obtain a series $K=K_0\supset K_1\supset \cdots\supset K_r$ of 
normal subgroups of $G$ and characters $\nu_i\in \Irr(K_i)$ such 
that $(K_i, K_{i+1}, \nu_i)$ is an extreme triple  for $0\leq i<r$ 
and where $\nu_r$ is linear. As before, we have that
$r< n$ and $K_r/\Ker(\chi\psi)$ is abelian.

 For $0\leq i<r$, let $L_i\triangleleft G$ such that $K_i\subset L_i\subseteq
K_{i_1}$ 
and $L_i/K_i$ is a chief factor of $G$. Since 
$G$ is supersolvable, we have that $L_i/K_i$ is cyclic of 
prime order and therefore $G/{\bf C}_G(L_i/K_i)$ is abelian.
Thus $\dl(K_{i}/K_{i+1})\leq 2$ for $0\leq i<r$. 
Therefore
$$\dl(K/\Ker(\chi\psi))\leq \sum_{i=1}^r \dl(K_{i-1}/K_{i})+
\dl(K_r/\Ker(\chi\psi))\leq 2r+1.$$
Since $r<n$, we conclude that 
$\dl(K/\Ker(\chi\psi) \leq 2(n-1) + 1= 2n-1$ and the proof is complete. 
\end{proof}
\end{section}
\begin{section}{Proof of Theorem B}
\begin{proof}[Proof of Theorem B]
 Set $K= \Ker(\alpha)$.
Let $\theta \in \Irr(K)$ be as in Lemma \ref{partb}. Since 
$[\chi_K, \theta]\neq 0$ and $[\overline{\psi_K}, \theta]\neq 0$, and
$K$ is normal in $G$, we have that 
$\theta(1)$ divides both $\chi(1)$ and $\psi(1)$. 
 Since  $(\chi(1), \psi(1))=1$, it follows that  
 $\theta(1)=1$.
 
Observe that $[K,K]\subseteq \Ker(\theta)\subseteq K$ since $\theta(1)=1$.
Also observe that $[K,K]$ is normal in
 $G$ since $K$ is normal in $G$. 
Since $\theta$ is an irreducible constituent of $\chi$, $[K,K]\subseteq
\Ker(\theta)$ and
$[K,K]$ is normal in  $G$, we have that $\chi_{[K,K]}=\chi(1)1_{[K,K]}$.
Similarly we can check that $\psi_{[K,K]}=\psi(1)1_{[K,K]}$.
Thus 
 $[K, K] \subseteq \Ker(\chi\psi) \subseteq K$.
Therefore $K/ \Ker(\chi \psi)$ is abelian and thus
 $\dl (\Ker(\alpha)/ \Ker(\chi\psi))\leq 1$. 
 \end{proof}

\begin{example} Let $G$ be a solvable group and $\chi$, $\psi \in \Irr(G)$.
Assume that $$\dl(\Ker(\alpha)/ \Ker(\chi\psi))\leq 1$$
\noindent  for all 
irreducible constituents $\alpha$ of the product $\chi \psi$. 
Then we do not necessarily have that $(\chi(1), \psi(1))=1$. 
\end{example}
\begin{proof} Let $p$ be an odd prime. 
Let $E$ be an extra-special group of order $p^3$. Choose
a character $\chi\in \Irr(E)$  such that $\chi(1)=p$. 
Let $\psi= \overline{\chi}$. Observe that $\Ker(\chi\psi)= {\bf Z}(E)$.
Also observe that $(\chi(1), \psi(1))= \chi(1)=p$. Since
$E/{\bf Z}(E)$ is abelian, for all the irreducible 
constituents $\alpha$ of the product $\chi\psi$,
we have that  $\dl(\Ker(\alpha)/ \Ker(\chi\psi))=1$.
\end{proof}
\end{section}
\begin{section}{Further results}

Assume that $G$ is a finite group and $\chi, \psi$ are 
irreducible characters of $G$
such that the product
$\chi\psi$ has a linear constituent $\alpha$. 
Observe that $ \alpha \overline{\chi} \in \Irr(G)$ since 
$\alpha (1)=1$ and $\chi \in \Irr(G)$.
Since $[\psi, \alpha \overline{\chi}]=[\chi\psi, \alpha]>0$, 
we have that $\psi = \alpha \overline{\chi}$. In particular
 $\chi(1)=\psi(1)$. If, in addition, the group $G$ is solvable, more
information is available about $\chi(1)$. 
Since $\alpha(1)=1$, by Exercise 4.12 of \cite{isaacs}
 we have that $[\chi\psi,\alpha]=1$.
If $\chi\psi(1)\neq 1$, by Proposition \ref{ai} there exists 
a character 
$\beta\in \Irr(G)$ such that $\beta\neq \alpha$ and $[\beta,\chi\psi]=1$.

The main results of this section are Propositions \ref{chi(1)factors} and 
\ref{ai}. Those are corollaries of Theorem B and C
of \cite{edith} and Lemma \ref{thesame}.

\begin{lemma}\label{thesame}
Let $G$ be a finite group and $\chi\psi \in \Irr(G)$. Set
\begin{equation*}\label{sumapsi}	
 \chi\psi=  \sum_{i=1}^n a_i \alpha_i
\end{equation*}
\noindent where $n>0$ and 
$a_i>0$ is the multiplicity of 
$\alpha_i\in \Irr(G)$ in $\chi\psi$, for each $i=1, \ldots, n$.
 
If $\alpha_1(1)=1$, then the irreducible constituents of 
the character $\chi\overline{\chi}$ are 
$1_G$, $ \overline{\alpha_1}\alpha_2$, $\overline{\alpha_1}\alpha_2, \ldots, 
\overline{\alpha_1}\alpha_n$, and 
\begin{equation*}
 \chi\overline{\chi}= a_1 1_G + \sum_{i=2}^n a_i (\overline{\alpha_1}\alpha_i)
\end{equation*}
\noindent where $n>0$ and 
$a_i>0$ is the multiplicity of 
$\overline{\alpha_1}\alpha_i$ in $\chi\overline{\chi}$. In particular, 
$\eta(\chi\psi)= \eta(\chi\overline{\chi})$.
\end{lemma} 
\begin{proof}
Observe that $ \alpha_1\overline{\chi} \in \Irr(G)$ since 
$\alpha_1(1)=1$ and $\chi \in \Irr(G)$.
Since $[\psi, \alpha_1 \overline{\chi}]=[\chi\psi, \alpha_1]>0$, 
we have that $\psi = \alpha_1 \overline{\chi}$.

Since $\psi=\alpha_1 \overline{\chi}$  and $\alpha_1$ is a linear character,
we have that $(\alpha_1)^{-1}=\overline{\alpha_1}$ and 
$$\chi\overline{\chi}= \sum_{i=1}^n a_i (\overline{\alpha_1} \alpha_i).$$
Observe that $\overline{\alpha_1} \alpha_i \in \Irr(G)$ since 
$\alpha_1$ is a linear character and $\alpha_i \in \Irr(G)$.
  Also observe that 
$\overline{\alpha_1} \alpha_i\neq \overline{\alpha_1} \alpha_j$ if
$\alpha_i \neq \alpha_j$. Thus the distinct irreducible 
constituents of $\chi\overline{\chi}$ are 
$1_G= \overline{\alpha_1}\alpha_1$, 
$ \overline{\alpha_1}\alpha_2$, $\overline{\alpha_1}\alpha_2, \ldots, 
\overline{\alpha_1}\alpha_n$, and $a_i$ is the multiplicity of 
$\overline{\alpha_1} \alpha_i$ in $\chi\overline{\chi}$.
 \end{proof}
\begin{proposition}\label{chi(1)factors}
Let $G$ be a finite solvable group. Let $\chi, \psi \in \Irr(G)$ 
with $\chi(1)>1$. If the product $\chi\psi$ has a linear 
constituent,  then $\chi(1)=\psi(1)$ and $\chi(1)$ has at most 
$\eta(\chi \psi)-1$ distinct
prime divisors. 

If, in addition, $G$ is supersolvable,
then $\chi(1)$ is a product of at most $\eta(\chi \psi)- 2$ primes.
\end{proposition}
\begin{proof}
By Theorem C of \cite{edith}, we have that $\chi(1)$ has at most 
$\eta(\chi\overline{\chi})-1$ distinct prime divisors. Also by Theorem C of
\cite{edith} if, in addition, 
$G$ is supersolvable,
then $\chi(1)$ is a product of at most $\eta(\chi\overline{\chi})- 2$ primes.
Thus the result follows from Lemma \ref{thesame}.
\end{proof}

\begin{proposition}\label{ai}
Let  $G$ be a finite solvable group and $\chi, \psi \in \Irr(G)$ 
with $\chi(1)>1$. Assume that the product $\chi\psi$ has a linear
constituent $\alpha_1$ and the  decomposition of
the character $\chi \psi$ into its distinct 
irreducible 
constituents $ \alpha_1, \ \alpha_2, \ldots, \alpha_n$
has the form
\begin{equation*}	
 \chi\psi=  \sum_{i=1}^n a_i \alpha_i
\end{equation*} 
\noindent where $n>0$ and 
$a_i>0$ is the multiplicity of 
$\alpha_i$ for each $i=1, \ldots, n$.
 
Then $a_1=[\alpha_1, \chi\psi]=1$ and 
$1 \in \{ a_i \mid i=2, \ldots, n\}$. 
\end{proposition}
 \begin{proof} 
Let $\{\theta_i \in \Irr(G)^{\#}\mid i=2, \ldots , \eta(\chi\overline{\chi})
\}$
be the set of non-principal
 irreducible constituents of $\chi \overline{\chi}$.
 If
 $\, \Ker(\theta_j)$ is maximal under inclusion  among the 
 $\Ker(\theta_i)$
for $i=2, \ldots, n$, then 
$ [ \chi \overline{\chi} , \theta_j] =1$ by Theorem C of \cite{edith}.
Thus $1\in \{ [\chi \overline{\chi} , \theta_i] \mid 
 i=2, \ldots, n\}$. The result then follows from Lemma \ref{thesame}. 
\end{proof}

\end{section} 
 
\bibliographystyle{amsplain}

\end{document}